\documentclass[12pt,leqno,twoside]{amsart}
\usepackage{amssymb}
\usepackage{latexsym}
\usepackage{verbatim}
\usepackage{epsfig}
\usepackage[all]{xy}

\newtheorem{theorem}{Theorem}[section]
\newtheorem{corollary}[theorem]{Corollary}
\newtheorem{lemma}[theorem]{Lemma}
\newtheorem{proposition}[theorem]{Proposition}
\theoremstyle{definition}
\newtheorem{definition}[theorem]{Definition}
\theoremstyle{remark}

\theoremstyle{remark}

\theoremstyle{remark}

\theoremstyle{remark}

\theoremstyle{remark}

\theoremstyle{remark}

\renewcommand{\Box}{\square}    
\renewcommand{\Bbb}{\mathbb}
\newcommand{\cal}{\mathcal}

\newcommand{\Cone}{{\rm{Cone}}}
\newcommand{\hd}{{\rm{hd}}} 
\newcommand{\codim}{{\rm{codim}}}  
  
\renewcommand{\rhd}{{\rm{rhd\hspace{1.5pt}}}}
\newcommand{\h}{{\rm{ht}}}

\renewcommand{\int}{{\rm{int}}}

\newcommand{\Sing}{{\rm{Sing}}}

\newcommand{\id}{{\rm{id}}}

\newcommand{\Iso}{{\rm{Iso}}}
\newcommand{\im}{\mathop{\rm{im}}\nolimits}

\newcommand{\cl}{{\rm{closure}}}

\newcommand{\e}{\varepsilon}
\newcommand{\fin}{\hspace*{\fill}$\Box$}
\newcommand{\m}{\setminus}

\newcommand{\hvar}{{\rm{hvar}}}

\newcommand{\cD}{{\cal D}}

\newcommand{\cH}{{\cal H}}
\newcommand{\cS}{{\cal S}}

\newcommand{\cT}{{\cal T}}

\newcommand{\cW}{{\cal W}}

\newcommand{\cN}{{\cal N}}
\newcommand{\cY}{{\cal Y}}

\newcommand{\bC}{{\Bbb C}}

\newcommand{\bP}{{\Bbb P}}

\newcommand{\bX}{{\Bbb X}}
\newcommand{\bY}{{\Bbb Y}}

\begin{document}

\title[Higher homotopy variation]{The Zariski-Lefschetz principle for higher homotopy
  groups of nongeneric pencils}
\author{Mihai Tib\u ar}

\date{\today}

\address{Math\' ematiques, UMR 8524 CNRS,
Universit\'e des Sciences et Technologies de Lille, \  59655 Villeneuve d'Ascq, France.}
\email{tibar@math.univ-lille1.fr}

 

\begin{abstract}
We prove a general
Zariski-van Kampen-Lefschetz type theorem for higher homotopy groups
of generic and  nongeneric pencils on singular open complex spaces.
\end{abstract}

\maketitle
\setcounter{section}{0}
\section{Introduction}


We find natural generators for
 the ``vanishing homotopy'' which occurs
 in a pencil of hypersurfaces on an analytic variety $X$. This result is the analogue
 of the Zariski-Lefschetz theorem on ``vanishing
cycles'' for homology groups. In the same time this generalises to
 higher homotopy groups the Zariski-van Kampen theorem
on the fundamental group. 

It turns out that our geometric method allows to prove such a result
 in a general setting, which includes for instance the use of nongeneric
 pencils, as we show in the following. Let  $Y$ be a compact complex analytic space
and let $V\subset Y$ be a closed analytic subspace
such that $X := Y\m V$ is connected, of pure complex dimension $n\ge 3$.
 A {\em pencil} on $X$ is the restriction to $X$ of a pencil on $Y$. 
 The latter is by definition the ratio of two 
sections, $f$ and $g$, of some holomorphic line bundle over $Y$. The
 meromorphic function $h:= 
f/g : Y \dashrightarrow \bP^1(\bC)$ is holomorphic over the 
complement $Y\m A$ of the axis $A := \{ f=g=0\}$.
  To define {\em singularities} of pencils in a generalised sense, we first consider a Nash 
blowing-up $\bY$ along the axis (see \S\ref{s:sing}), 
with projections $p \colon \bY \to \bP^1$ and
$\sigma \colon \bY \to Y$.  Then there is a well-defined singular locus of $p$ with
respect to a certain natural Whitney stratification 
$\cS$, denoted by $\Sing_\cS p$.

We shall assume throughout this paper that $\dim \Sing_\cS p \le 0$, i.e.
 the singularities are {\em isolated}. It should be pointed out that
 singularities may lie on the axis, i.e. that one may have
 $A\cap \sigma(\Sing_\cS p)\not= \emptyset$.
 This
 occurs for instance if the transversality of the 
axis to strata of $Y$ fails at finitely many points: see
Proposition \ref{p:ijm} and  \cite[\S 2]{Ti-lef}, \cite[\S 2]{Ti-xiv} for
more details.

Such nongeneric pencils were introduced in our preprint
\cite{Ti-pre} and this viewpoint allowed us 
to prove connectivity theorems of Lefschetz type \cite{Ti-lef}, and 
a far reaching extension of the Second Lefschetz Hyperplane
Theorem \cite{Ti-xiv}.
 In the literature, examples of nongeneric pencils occured sporadically; 
more recently they
got into light, e.g. \cite{KPS, Ok2}, precisely because
 one can use nongeneric pencils towards more
efficient computations.
 For instance, M. Oka uses special pencils tangent to flex
points of projective curves in order to compute the fundamental group of the
complement, see e.g. \cite{Ok1, Ok2}. 
Nongenericity also proved to be a useful concept for treating the topology of
 polynomial functions (see \cite{NN2, Ti-xiv}) and of complements of arrangements
(e.g. \cite{LT}). A highly nongeneric situation occurs when $V$ 
contains a member of the pencil, see \S \ref{ss:poly}.

In order to express how the ``homotopy vanishing cycles'' are generated we need
to introduce a {\em homotopy variation map}. Roughly, this goes as follows
(see \S\ref{ss:var}).
Let $a\in\Lambda = p(\Sing_\cS p)$ be  some critical value of the pencil and
denote by $X_a := X \cap p^{-1}(a)$ the corresponding critical fibre. 
For a general fibre of the pencil $X_c$
 near to $X_a$, we denote 
by $X^*_{a}$ the space $X_{c}\m \cup_j B_{j}$, where $B_{j}$ is a Milnor ball
at some point of the set $X_a\cap\sigma(\Sing_\cS p)$.
We construct in \S \ref{ss:var} a global variation map
for higher homotopy groups:
\begin{equation}\label{eq:hvar} 
\hvar_a : \pi_q(X_{c}, X^*_{a}, \centerdot) \to \pi_q(X_{c}, \centerdot),
\end{equation}
 for any $q\ge 3$. Under the condition that the involved spaces are 
path-connected and that
$\pi_1(X_{c}, X^*_{a}, \centerdot)$ is trivial, this map commutes with the
actions of $\pi_1(X^*_{a}, \centerdot)$ and of $\pi_1(X_{c},
\centerdot)$ via the surjection $\pi_1(X^*_{a}, \centerdot) \twoheadrightarrow
\pi_1(X_{c},\centerdot)$, respectively. 
 
In case of homology groups, one encounters variation maps in the literature ever
 since Picard and Lefschetz 
studied the monodromy around an isolated singularity of a 2-variables 
holomorphic function \cite{PS, Lef}. Generalised variation maps and
 Picard-Lefschetz formulas play a key role
 in studying the topology of pencils of hypersurfaces 
(e.g. \cite{Ph, Mi, La, Lo, Ga, NN2},   \cite{Va1, Va2}). 
We have introduced in \cite{Ti-xiv} a global variation map
in homology which generalises the local variation map to pencils with isolated
singularities on arbitrarily singular (open) underlying
spaces.

 Let us now state our general Zariski-van Kampen-Lefschetz
 type result which uses the homotopy variation map (\ref{eq:hvar}). 
We refer to
 \S \ref{s:sing} and  \S \ref{ss:var} for the precise definitions
 of all the ingredients. The homotopy depth condition ``$\hd$'' is 
recalled in Definition \ref{d:hd} and the notation ``$X_{a}^*$'' is explained
 in \ref{ss:notation}.
  

\begin{theorem}\label{t:main}
  Let $h \colon Y \dashrightarrow \bP^1$ define a pencil with isolated
  singularities such that the axis A of the pencil is not included in $V$ and
  that the general fibre $X_c$ on $X = Y\m V$ is path
  connected.  
   Let us consider the following three conditions:
\begin{itemize} 
\item[(i)] $(X_c, X_c\cap A)$  is $k$-connected;
\item[(ii)] $(X_c, X_{a}^*)$  is $k$-connected, for any critical value
  $a$ of the pencil; 
\item[(iii)]  $\hd_{\bX \cap \Sing_\cS p} \bX \ge k+2$.  
\end{itemize} 
  Then we have the following two conclusions:
\begin{enumerate}
\item If for some $k\ge 0$ the conditions \rm (i)\it, \rm (ii)\it\  and \rm (iii)\it\ hold then 
$\pi_q(X,X_c, \centerdot) =0$ for all $q\le k+1$. 
\item If for some $k\ge 2$ the conditions \rm (i)\it\ and \rm (ii)\it\ hold,  and if condition
  \rm (iii)\it\ is  replaced by the following one:
\newline \hspace*{-5mm}
 \rm (iii')\it\ $\hd_{\bX \cap \Sing_\cS p} \bX \ge k+3$, \\
then the kernel of the surjective map $\pi_{k+1} (X_c, \centerdot) \twoheadrightarrow \pi_{k+1}(X,
\centerdot)$ is 
generated by the images of the variation maps $\hvar_a$. 
\end{enumerate}
\end{theorem}
The conditions in our theorem are satisfied by a large class of spaces $X$.
For instance, the conditions (ii) and (iii), resp. (iii'), are both fulfilled as soon as
$X$ is a singular space which is a locally complete intersection of dimension
$n\ge k+2$, resp. $n\ge k+3$ (see \S \ref{ss:depth}). We send to
\S\ref{s:further} for many other comments on these
conditions and on what they become in more particular situations.

 In what concerns Part (a) of Theorem \ref{t:main}, this is a general
 ``first Lefschetz hyperplane theorem'' for homotopy groups. It has  been proved
 in some particular cases in the past and it
recovers, with significantly
 weakened hypotheses, our recent result \cite{Ti-lef}. 

The focus of this paper
 is actually on  Part (b) of Theorem \ref{t:main}, which treats the 
``homotopy vanishing cycles''. 
 Our approach is in the spirit of the Lefschetz method 
\cite{Lef}, as presented by Andreotti and Frankel in \cite{AF}.
Since the homology conterpart of our result (which is easier) 
is presented in \cite{Ti-xiv},
 we concentrate here on the more delicate constructions and ingredients of the
 proof which are needed in case of homotopy groups.
We use homotopy excision (the Blackers-Massey theorem), more precisely
 the general version proved by Gray \cite[Cor; 16.27]{Gr}, 
and the comparison between homotopy and
homology groups via the Hurewicz map. Our result does not cover the case
 $k=1$, hence the homotopy group $\pi_2$, because of the technical restriction
 in the definition of the variation map (\ref{eq:hvar}).

The cornerstone of the proof is perhaps Proposition \ref{p:4} in which we identify the
 homotopy variation map to the boundary map of the long exact sequence of the pair
$(\bX^*_D, X_c)$, by expanding Milnor's explanation of the Wang exact sequence  
\cite[pag. 67]{Mi}. All those techniques are coupled with the general theory
of stratified singular spaces and of Milnor fibrations of stratifies
singularities of holomorphic functions (see e.g. \cite{Le, GM}). Notably
Proposition \ref{p:3} is a key interplay between the use of the singular  Milnor
fibration and the homotopy excision.\\

 Theorem \ref{t:main} represents a synthetic viewpoint on Zariski-van Kampen
type results for higher homotopy groups, having in the background the
pioneering work of  
Lefschetz \cite{Lef}, Zariski and van Kampen \cite{vK}. It recovers several
 cases considered before
 in the literature and in the same
 time it extends the range of applicability of the
 Zariski-van Kampen-Lefschetz principle.
For instance we show in \S \ref{ss:classic} how to get a version of Theorem
\ref{t:main} in the case of complements of singular projective
spaces, which seems to be a new result.
  Some other possible applications of Theorem \ref{t:main} 
 parallel the ones for homology
 groups which have been discussed 
in \cite[Applications]{Ti-lef} and \cite{LT}.

 Such a progress would not have been possible without the contribution
of a long list of articles around the Lefschetz slicing principle,
from which we took a part of the inspiration, such as \cite{AF, Ch-ens, Ch, Fu,
  FL, GM, HL-lef, HL-gen, HL, La, Lef, Li, Lo, Mi}.
Results of Zariski-van 
Kampen type have been proved by Libgober for 
generic Lefschetz pencils of hyperplanes and in the particular
case when $X$ is the 
complement in $\bP^n$ or $\bC^n$ of a hypersurface with isolated singularities  
\cite{Li}.  D. Ch\'eniot, A. Libgober and C. Eyral  authored in 2002 two preprints
 \cite{CL, CE} treating, with different background, interesting particular
 cases of the problem we are concerned with here.

This paper completes in detail a certain  part of the
  program ``Lefschetz principle for nongeneric pencils'' presented in our 2001
 preprint \cite{Ti-pre}.
 We own warm thanks to several institutions from which we have got
 hospitality in connection to this project: Newton Institute (2000), IAS
 Princeton (2002), CRM Barcelona (2004).
\section{Pencils with singularities in the axis}\label{s:sing}

As denoted before, let $X := Y\m V$ be the difference of two compact complex
 analytic spaces.
We introduce the basic notations and recall from \cite{Ti-lef} what are
 the singularities of a pencil, in a stratified sense. 

Consider the compactification of the graph of $h$, namely the space:
\[ \bY :=\cl \{ y \in Y, \ [s;t] \in \bP^1 \mid sf(y) - t g(y)  =0\} 
\subset Y\times \bP^1. \]
Let us denote $\bX:= \bY \cap (X\times \bP^1)$. Consider the projection $p : \bY \to 
\bP^1$,
its restriction $p_{|\bX} : \bX \to \bP^1$ and
the projection to the first factor $\sigma : \bY \to Y$.
 Notice that the restriction of $p$ to $\bY\m(A\times \bP^1)$ can be identified with $h$.

Let $\cW$ be an analytic Whitney stratification of $Y$ such that $V$ is union of strata. 
Its restriction to the open set $Y\m A$ induces a Whitney 
stratification on $\bY\m(A\times \bP^1)$, via the mentioned identification.
We denote by $\cS$ the coarsest analytic Whitney
stratification on $\bY$ which coincides over $\bY\m (A\times\bP^1)$ with
the one induced by $\cW$ on $Y\m A$. This stratification exists within a neighbourhood of 
$A\times \bP^1$, by usual 
arguments (see e.g. \cite{GLPW}), hence such stratification is well defined on $\bY$.  We 
call it the
{\em canonical stratification} of $\bY$ generated by the stratification $\cW$ 
of
$Y$. 


\begin{definition}\label{d:sing}\cite{Ti-lef}
 One calls {\em singular locus} of $p$ with respect to $\cS$ the following subset of 
$\bY$:
 \[  \Sing_\cS p := \bigcup_{\cS_\beta \in \cS} \Sing
p_{|\cS_\beta}.\]
The set $\Lambda := 
p(\Sing_\cS p)$ is called the set of {\em critical values} of $p$ with respect to $\cS$.
We say that the pencil defined by $h:= 
f/g : Y \dashrightarrow \bP^1$ is a {\em pencil 
with isolated singularities} if $\dim \Sing_\cS p \le 0$.
\end{definition}

Remark that $\Sing_\cS p$ is  closed analytic.  Remark also that $p$ is proper 
and analytic in a neighbourhood of $A\times \bP^1$ and that $\cS$ has finitely 
many strata. Then, by a Bertini-Sard result which
is a consequence of the First Isotopy Lemma \cite{Th}, see also \cite{Ve},
 it follows that $\Lambda \subset \bP^1$ is a finite set and that the maps  $p: 
\bY\m p^{-1}(\Lambda) \to 
\bP^1\m \Lambda$ and $p_{|\bX}:
\bX \m p^{-1}(\Lambda) \to \bP^1\m \Lambda$ are  stratified
locally trivial fibrations. In particular, $h: Y\m (A\cup h^{-1}(\Lambda)) \to \bP^1\m 
\Lambda$ is a 
locally trivial fibration.

   
    Such nongeneric pencils may occur whenever the axis $A$ is not in general
    position at finitely many points. More precisely, we have the
    following example of situation:

\begin{proposition}\label{p:ijm} {\rm \cite[Prop. 2.4]{Ti-lef}}\\
Let $X= Y\m V \subset \bP^N$ be a quasi-projective (singular) variety and let
$\hat h= \hat f/\hat g$ define a pencil of hypersurfaces in $\bP^N$ with axis $\hat A$.  Let
$C$ denote the set of points on $\hat A \cap Y$ where some member of the
pencil is singular or where $\hat A$ is not stratified transversal to $Y$. If
$\dim C\le 0$ and if the stratified critical points of the restriction $\hat
h_{|} \colon Y\m A\to \bP^1$ are isolated, then $\dim \Sing_\cS p \le 0$.
\fin
\end{proposition}

 We shall assume throughout the paper that our pencil on $X$ has isolated
 singularities. Then, by the compactness of $Y$, the set $\Sing_\cS p$ is a
 finite set.  Let us emphasize that the points of $\Sing_\cS p$ are not
 necessarily contained in $\bX$ and that some of them may be
 on the blown-up axis $A\times \bP^1$.

\section{Homotopy variation maps}\label{ss:var}

 Inspired by the homology constructions throughout the literature
 and by a particular Zariski-van-Kampen type theorem for higher homotopy
 groups proved in the beginning of the 1990's by  Ligbober \cite{Li}, 
we have given in the preprint \cite{Ti-pre} (unpublished) the main lines
 of the construction of {\em global homotopy
variation maps} for pencils with isolated singularities. 
We complete here in full detail this construction\footnote{A similar construction in case of
 homology groups appears in  \cite{Ti-xiv}.}.

 For any $M\subset \bP^1$, we denote $\bY_M := p^{-1}(M)$ and $\bX_M 
:= \bX \cap \bY_M$ and $X_M := \sigma(\bX_M)$.
Let $a_i \in\Lambda$ and let  $\Sing_\cS p = \{a_{ij}\}_{i,j} \subset \bY$, where
$a_{ij}\in \bY_{a_i}$.  For $c\in \bP^1\m \Lambda$ we say that $\bY_c$,
respectively $\bX_c$, respectively 
$X_c$, is a general fibre.

\subsection{Milnor fibration at isolated stratified critical points of
  holomorphic functions}\label{ss:milnor} 

At some singularity $a_{ij}$, in local coordinates, we take a small ball
 $B_{ij}$ centered 
at $a_{ij}$. For small enough radius of $B_{ij}$,
this is a {\em Milnor ball}\footnote{One may consult \cite{GM, Le}
for the definition and the terminology concerning the local Milnor fibration in case
of a holomorphic function on a singular analytic space, at an isolated 
critical point with respect to a given stratification of the space.}  of the
 local holomorphic germ of the function $p$ at
 $a_{ij}$.  Next 
we may take a small enough disc $D_i\subset \bP^1$ at $a_i \in \bP^1$, 
so that $(B_{ij}, D_i)$ 
is {\em Milnor data} for $p$ at $a_{ij}$. Moreover, we may do this for all (finitely many) 
singularities in the fibre $\bY_{a_i}$, keeping the same disc $D_i$, provided
 that it is small enough.

Since the function $p : \bY \to \bP^1$ has isolated stratified singularities,
the fibres of $p$ are endowed with the stratification induced by $\cS$ except
that we have to introduce the point-strata $\{a_{ij}\}$. Every fibre has
therefore a natural induced Whitney stratification.

By the general theory of holomorphic function germs \cite{GM, Le} we have that
the restriction $p_| : \bY_{\partial \bar D_i} \to \partial \bar D_i$
 is a  {\em locally trivial} stratified fibration over the circle $\partial \bar D_i$.
Since $p : \bY \to \bP^1$ has a stratified 
 isolated singularity at $a_{ij}$, the sphere $\partial \bar
 B_{ij}$ is stratified transversal to all the fibres of $p$ above the points of
 $D_i$. It follows that the restriction of $p$ to the pair of spaces 
 $(\bY_{D_i} \m \cup_j B_{ij}, \bY_{D_i} \cap \cup_j \partial \bar B_{ij})$ is
 a {\em trivial stratified fibration} over $D_i$.

\subsection{Geometric monodromy}\label{ss:geomono} 
Let's choose $c_i \in \partial \bar D_i$. 
One constructs in the usual way (see Looijenga's
  similar discussion in \cite[2.C, page 31]{Lo}), by using a stratified vector field,
 a characteristic morphism $h_i : X_{c_i} \to
X_{c_i}$ of the fibration over the circle $p_| : \bX_{\partial \bar D_i} \to
\partial \bar D_i$, 
which is a stratified homeomorphism and it is the identity on $X_{c_i}\m 
\cup_j B_{ij}$. This is called the {\em geometric monodromy} of the locally
trivial fibration
$p_| : \bX_{\partial \bar D_i} \to \partial \bar D_i$ and corresponds to one 
counterclockwise loop around the circle $\partial \bar D_i$.
  Clearly $h_i$ is not uniquely defined with these properties, but its
  relative isotopy class in the group $\Iso (X_{c_i}, X_{c_i}\m \cup_j
  B_{ij})$ of relative isotopy classes of stratified homeomorphisms which are
  the identity on  $X_{c_i}\m \cup_j B_{ij}$ is unique.
 
 \subsection{Construction of the homotopy variation map}\label{ss:hvar} 
 We assume in the rest of this section 
that $X_{c_i}\m \cup_j B_{ij}$ is path connected, for all $i$. This assumption
is actually fulfilled under the hypotheses of our Theorem \ref{t:main}(b). More precisely,
the fact that the general fibre $X_c$ is path connected and that
 $(X_{c_i}, X_{c_i}\m \cup_j B_{ij})$ is
1-connected imply that $X_{c_i}\m \cup_j B_{ij}$ is path connected.
Let $\gamma \colon (D^q, S^{q-1}, u) \to (X_{c_i}, X_{c_i}\m \cup_j
 B_{ij}, v)$, for $q\ge 3$,  be some 
continuous map, and let $[\gamma]$ be its homotopy class in $\pi_q(X_{c_i}, 
X_{c_i}\m \cup_j B_{ij}, v)$. 
Note that $[h_i\circ \gamma]$ is a well-defined element of $\pi_q(X_{c_i}, 
X_{c_i}\m \cup_j B_{ij}, v)$ and does not depend on the representative $h_i$
in its  relative isotopy class. Consider the map $\rho \colon (D^q, S^{q-1},
 u) \to (D^q, S^{q-1}, u)$ which is 
the reflection into some fixed generic hyperplane through the origin, which 
contains $u$. Then $[\gamma 
\circ \rho]$ is the inverse  $-[\gamma]$ of $[\gamma]$. We use the additive
notation since the relative homotopy groups $\pi_q$ are abelian for $q\ge 3$.
 
 Consider now the map:
\begin{equation} \label{eq:mu}
\mu_i \colon  \pi_q(X_{c_i}, X_{c_i}\m \cup_j 
B_{ij},v) \to 
\pi_q(X_{c_i}, X_{c_i}\m \cup_j B_{ij},v)
\end{equation}
defined as follows:
 $\mu_i([\gamma]) = [(h_i\circ \gamma) *(\gamma \circ \rho)] = [h_i\circ
 \gamma] - [\gamma]$.

The notation ``$*$'' stands for the operation on maps which induces the
group structure of the relative $\pi_q$. The map $\mu_i$ is well-defined and
it is an automorphism of $\pi_q(X_{c_i}, X_{c_i}\m \cup_j 
B_{ij},v)$, because of the abelianity of the relative homotopy groups
 $\pi_q$, for $q\ge 3$. For the time being, we cannot extend this construction
 to $q=2$. In case $q=1$ there are the monodromy relations which enter in the well-known
 Zariski-van Kampen theorem. 

\smallskip
\noindent
{\bf Claim}: \it the map $\nu_i := (h_i \circ \gamma) * (\gamma\circ \rho)
 \colon (D^q,S^{q-1},u) \to
(X_{c_i}, X_{c_i}\m \cup_j B_{ij},v)$ is homotopic, relative to $S^{q-1}$,
 to a map $(S^q, u)\to (X_{c_i}, v)$. \rm \smallskip

To prove the claim, we first observe that for the restriction of $\gamma$ to
 $S^{q-1}$, which we shall call 
$\gamma'$, we have $h_i \circ \gamma' = \gamma'$, since the geometric
 monodromy $h_i$ is the identity on $\im \gamma'$. Without loss of generality,
we may and shall suppose in the following that the geometric monodromy $h_i$ is the
 identity on a small tubular neighbourhood $\cT$ of $X_{c_i}\m \cup_j B_{ij}$ 
within $X_{c_i}$.

 Let $\nu_i'$ denote the restriction of $\nu_i$ to $S^{q-1}$. We claim that
 one can shrink $\im \nu_i'$ to the point $v$ through a homotopy within $\im
 \nu_i'$. The precise reason is that, by
 the definition of $\rho$, the map  $\gamma' * (\gamma' \circ \rho)$ is
 homotopic to the constant map to the base-point $v$, 
through a deformation of the image of $\gamma' * (\gamma' \circ \rho)$ within itself.
 We may moreover extend this homotopy to a 
continuous deformation of the map $\nu_i$
 in a small collar neighbourhood $\tau$ of $S^{q-1}$ in $D^q$, which has its
 image in the tubular neighbourhood $\cT$, such that it is the identity on the interior
boundary of $\tau$ (which is a smaller sphere $S_{1-\e}^{q-1}\subset D^q$).
 Then we extend this deformation as the identity to the
 rest of $D^q$. 
In this way we have constructed a homotopy between $\nu_i$ and a map 
$(D^q,S^{q-1},u) \to (X_{c_i},v,v)$, which is nothing else than a map $(D^q,u) \to
(X_{c_i},v)$. This represents an element of the group $\pi_q(X_{c_i},v)$, 
so our claim is proved.

This construction is well-defined up to homotopy equivalences, hence it yields
a map:
\begin{equation}\label{eq:vvar}
 \hvar_i : \pi_q(X_{c_i}, X_{c_i}\m \cup_j B_{ij},v) \to \pi_q(X_{c_i},v),
\end{equation} 
which we call the {\em (homotopy) variation map} of $h_i$. 

Let us point out  that $\hvar_i$ is
 a morphism of groups. Indeed, due to the abelianity we have:
  \[ \hvar_i([\gamma] +[\delta]) = [(h_i \circ
 \gamma) * (h_i\circ \delta)] - [\gamma * \delta] = [h_i \circ\gamma] +
 [h_i\circ \delta] -[\gamma] - [\delta ] \]
\[ = [h_i \circ\gamma]-[\gamma] +
 [h_i\circ \delta] -[\delta ] = \hvar_i([\gamma]) + \hvar_i([\delta]).\] 

By its geometric definition, this variation map enters in the following
commuting diagram:
\begin{equation}\label{eq:diag}
 \begin{array}{ccc}
 \pi_q(X_{c_i}, X_{c_i}\m \cup_j B_{ij},v) & \stackrel{\mu_i}{\longrightarrow}
 & \pi_q(X_{c_i}, 
X_{c_i}\m \cup_j B_{ij},v) \\
   \ \ \ \ \ \ \ \ \ \ \ \ \ \ \ \  \mbox{\tiny $\hvar_i$} \searrow  & \ &
 \nearrow \mbox{\tiny 
$j_\sharp$} \ \ \ \ \ \ \ \ \ \ \ \ \ \ \ \ \ \ \ \ \ \\
  \ & \pi_q (X_{c_i},v) & \ 
 \end{array} 
\end{equation}
where $j \colon (X_{c_i},v) \to (X_{c_i},X_{c_i}\m \cup_j B_{ij},v)$ is the 
inclusion.

 \subsection{Convention of notations}\label{ss:notation}
 In local coordinates at $a_{ij}$, $\bY_{a_i}$ is a germ of a stratified
 complex analytic space; hence, for a small enough ball $B_{ij}$, the set
 $\bar B_{ij}\cap X_{a_i}\m \cup_j 
a_{ij}$ retracts to $\partial \bar B_{ij}\cap X_{a_i}$, by the local conical 
structure of stratified
analytic sets \cite{BV}.
  It follows that $X_{a_i}\m \cup_j a_{ij}$ is homotopy equivalent, by retraction, to 
 $X_{a_i}\m \cup_j B_{ij}$. We have seen in \S \ref{ss:milnor} that 
the restriction $p_| : \bX_{D_i} \m \cup_j B_{ij} \to D_i$ is a trivial 
stratified fibration. This implies that the restriction to $\bX$, i.e. 
 $p_| : \bX_{D_i} \m \cup_j B_{ij} \to D_i$ is a trivial 
stratified fibration too. Since $X_{a_i}\m \cup_j B_{ij}$ and $X_{c_i}\m
 \cup_j B_{ij}$ are both fibres of this fibration, we shall identify them via
 a fixed trivialisation along a chosen path in $D_i$. 

In the remainder of this paper, for the sake of simplicity, we shall denote 
the spaces  $X_{c_i}\m \cup_j B_{ij}$ and $X_{a_i}\m \cup_j B_{ij}$ by the 
same symbol $X^*_{a_i}$.
In particular,
 $X_{a_i}\m \cup_j B_{ij}$ becomes in this way a subspace of $X_{c_i}$. 

 We shall use in the remainder the following notation 
 for the variation map, only as a symbolic replacement for (\ref{eq:vvar}), 
but not having an intrinsic meaning:
\begin{equation}\label{eq:var} 
 \hvar_i : \pi_q(X_{c_i}, X^*_{a_i}, \centerdot) \to \pi_q(X_{c_i}, \centerdot),
 \end{equation}
for any $q\ge 3$. This is actually the notation used in the statement of
Theorem \ref{t:main}.

 \subsection{The action of $\pi_1$}\label{ss:action}
We have the standard action of $\pi_1(X_c, \centerdot)$ 
on the  higher homotopy groups $\pi_q(X_c, \centerdot)$ 
and on any relative  homotopy group where $X_c$ appears at the second
position. Moreover, the hypotheses of Theorem \ref{t:main}(b)
imply that  $X^*_{a_i}$ is path connected for all
$i$, as explained in \S \ref{ss:hvar}, and that $\pi_1(X^*_{a_i}, \centerdot) 
\stackrel{i_*}{\to} \pi_1(X_{c_i}, \centerdot)$ is an isomorphism.
We shall therefore identify the two actions modulo the isomorphism $i_*$.

  So one has the action of
 $\pi_1(X^*_{a_i}, \centerdot)$
 on those higher relative homotopy groups of pairs
where $X^*_{a_i}$ is the second space of the pair, and one also has the action 
of $\pi_1(X^*_{a_i}, \centerdot)$, via the isomorphism $i_*$, on 
$\pi_q(X_{c_i}, \centerdot)$ and on the higher relative homotopy groups of pairs
where $X_{c_i}$ is the second space of the pair.
 
Let's then prove that the image of $\hvar_i$ is invariant under the action
 of $\pi_1(X_{c_i}, \centerdot)$. For some  $\beta \in \pi_1(X^*_{a_i}, \centerdot)$, we have 
$\hvar_i(\beta [\gamma])= \beta [h_i\circ \gamma] - \beta[\gamma] = 
 \beta \hvar_i([\gamma])$
since,  by the functoriality, the action of $\beta$ commutes with the 
morphism induced by the geometric monodromy map
$h_i : X_{c_i}\to X_{c_i}$.

We shall denote in the following by $\hat \pi_q$ the quotient of the
 corresponding higher homotopy
group by the action of $\pi_1(X_c, \centerdot)$ or of  $\pi_1(X^*_{a_i},
 \centerdot)$, eventually through the isomorphism $i_*$.  
We shall use the same notation $\hat \nu$ for
 the passage to the quotients of some morphism of groups $\nu$ which commutes
 with the action of $\pi_1$.

\section{Proof of Theorem \ref{t:main}}\label{proof}
We shall prove in parallel (a) and (b). The only proof of (a) 
would be much simpler and can be easily extracted.

 Relative to the claim (b), let us first remark that the kernel of 
$\pi_{k+1} (X_c, \centerdot) \twoheadrightarrow \pi_{k+1}(X,
\centerdot)$ is invariant under the the $\pi_1(X_c, \centerdot)$-action on 
$\pi_{k+1} (X_c, \centerdot)$, by the functoriality of this action. 
So in order to prove (b) it suffices to show that the image of 
 $\hat \partial : \hat \pi_{k+2}(X, X_c, \centerdot) \to \hat 
\pi_{k+1} (X_c, \centerdot)$ is generated by the images 
$\im(\widehat{\hvar_i})$, since $\im(\hvar_i)$ is also invariant under
the $\pi_1$-action, see \S\ref{ss:action}. This is what we shall do in the last part of the
proof below.

When dealing with homotopy groups, we shall need to apply the homotopy
excision theorem of 
Blackers and Massey in the form given by Gray \cite[Cor; 16.27]{Gr}. 
 For the connectivity claim (a),
 we follow in the beginning the lines of \cite{Ti-lef}.

\smallskip

 Let  $A':= A\cap X_c$ 
and assume that $A' \not= \emptyset$. 
 Let $K\subset \bP^1$
be a closed disc with $K\cap \Lambda = \emptyset$ and let $\cD$ denote the
closure of its 
complement in $\bP^1$. We denote by $S:= K\cap \cD$ the common boundary, which
is a circle, 
and take a point $c\in S$.
Then take standard paths $\psi_i$ (non self-intersecting, non mutually
intersecting) from 
$c$ to $c_i$, such that $\psi_i \subset \cD \m \cup_iD_i$. The configuration 
$\cup_i(\bar 
D_i \cup 
\psi_i)$ is a deformation retract of $\cD$.
We shall identify the fibre $X_{c_i}$ to the fibre $X_c$,  by parallel transport 
along the path $\psi_i$.

We have the natural commutative triangle:
\begin{equation}\label{eq:2}
\begin{array}{ccc}
 \pi_{k+2} (X_\cD, X_c, \centerdot) & \stackrel{\mbox{\tiny $\iota_\sharp$}}
{\longrightarrow} & \pi_{k+2} 
(X, X_c, \centerdot) \\
 \ \ \ \ \  \ \ \ \ \ \ \mbox{\tiny{$\partial_0$}} \searrow & \  & \swarrow  
\mbox{\tiny{$\partial$}} \ \ \ \ \  \ \ \ \ \ \ \\
  \  & \pi_{k+1} (X_c, \centerdot) & \  \\
\end{array}
 \end{equation}
 where $\partial$ and $\partial_0$ are boundary morphisms and $\iota$ is
 the inclusion of 
pairs $(X_\cD, X_c)\hookrightarrow (X, X_c)$. 

One sais that an inclusion of pairs of topological spaces 
 $(N, N')\hookrightarrow (M, M')$ is a {\bf $q$-equivalence} if this inclusion induces 
an isomorphism 
of relative homotopy groups for $j<q$ and a surjection for $j=q$.
\begin{lemma}\label{l:im}
Let $(X_c, A')$ be $k$-connected for $k\ge 0$.  Then the inclusion of pairs
$\iota : (X_\cD, X_c) \hookrightarrow (X, X_c)$ is a $(k+2)$-equivalence, for
$k\ge 0$. 
In particular,  after taking 
$\pi_1(X_c, \centerdot)$-quotients in diagram (\ref{eq:2}),  we have:  $\im \hat \partial = 
\im \hat \partial_0$.
\end{lemma}
\begin{proof}
 Under
the assumption  $(X_c, A')$ is $k$-connected, $k\ge 0$, it follows from 
\cite[Proposition 3.2(i)]{Ti-lef} that $(X_S, X_c)$ is $(k+1)$-connected.
 From the exact sequence of the triple $(X_\cD, X_S, X_c)$, we deduce that 
the inclusion 
$(X_\cD, X_c)\hookrightarrow (X_\cD, X_S)$ is a $(k+2)$-equivalence.
 Next, \cite[Proposition 3.2(ii)]{Ti-lef} shows that, by homotopy excision,
 the inclusion of pairs $(X_\cD, X_S)\hookrightarrow (X, X_K)$
is a $(k+2)$-equivalence. Since $X_K$ retracts in a fiberwise way to $X_c$, the pair
$(X, X_K)$ is homotopy equivalent to $(X, X_c)$.
\end{proof}
We further go to the blow-up $\bX$ and get:
\begin{lemma}\label{l:1}
 The inclusion $(\bX_\cD, X_c) \to 
(X_\cD, X_c)$ is a homotopy equivalence.
\end{lemma}
\begin{proof}
 It is a simple fact that $X_\cD$ is homotopy equivalent 
to $\bX_\cD$ to which 
one attaches along $A'\times \cD$ the product $A'\times \Cone (\cD)$ 
(see e.g. \cite[Lemma 3.1]{Ti-lef}). Since $\cD$ is 
contractible, our claim follows.
\end{proof}

We also need the following result from \cite{Ti-lef}, which uses 
\cite[Lemma 3.1]{Ti-lef} and Switzer's result \cite[6.13]{Sw}:
\begin{lemma}\label{l:prop}{\rm \cite[Prop. 3.2(b)]{Ti-lef}}\\
If $(\bX_{D_i}, X_{c_i})$ is 
$(k+1)$-connected, $k\ge 0$,  for all $i$, 
then $(\bX_\cD, X_c)$ is $(k+1)$-connected.

\fin
\end{lemma}

For all $i$, let $\iota_i : (\bX_{D_i}, X_{c_i}) \to 
(\bX_\cD, X_c)$ denote the inclusion of pairs and let  
$\partial_i : \pi_{k+2} (\bX_{D_i}, X_{c_i}, \centerdot) \to  \pi_{k+1}
(X_{c_i}, \centerdot)$ denote the boundary 
morphism in the pair $(\bX_{D_i}, X_{c_i})$.
 Consider now the following commutative diagram, for $k\ge 1$:
 \begin{equation}
\begin{array}{ccc}
 \oplus_i \pi_{k+2} (\bX_{D_i}, X_{c_i}, \centerdot) & \  & \  \\
 \mbox{\tiny{$\sum {\iota_i}_\sharp$}} \downarrow &  \ &  \searrow 
\mbox{\tiny{$\sum \partial_i$}} \ \ \ \ \ \ \ \ \ \ \  \ \ \ \  \ \ \ \ \  \ \ \ \ \ \\
  \pi_{k+2} (\bX_{\cD}, X_c, \centerdot) & 
\stackrel{\mbox{\tiny{$\partial_0$}}}{\longrightarrow} &
 \pi_{k+1} (X_c,  \centerdot),
 \end{array}
 \end{equation}
 where, for all $i$, we identify $(X_{c_i}, \centerdot)$ to  $(X_{c}, \centerdot)$ together
 with the base points $\centerdot$
  by using the paths $\psi_i$,
as explained in the beginning of this section.  We use the additive notations 
$\sum \partial_i$ and $\sum {\iota_i}_\sharp$ since for $k\ge 1$ the groups
$\pi_{k+1} (X_c,  \centerdot)$ and $\pi_{k+2} (\bX_{\cD}, X_c, \centerdot)$ are
   abelian.

By the functoriality of the action of $\pi_1(X_c,\centerdot)$, we get the induced ``hat''
 diagram:
\begin{equation}\label{eq:3}
\begin{array}{ccc}
 \oplus_i \hat \pi_{k+2} (\bX_{D_i}, X_{c_i}, \centerdot) & \  & \  \\
 \mbox{\tiny{$\sum {\hat{\iota_i}}_\sharp$}} \downarrow &  \ &  \searrow 
\mbox{\tiny{$\sum \hat \partial_i$}} \ \ \ \ \ \ \ \ \ \ \  \ \ \ \  \ \ \ \ \  \ \ \ \ \ \\
 \hat  \pi_{k+2} (\bX_{\cD}, X_c, \centerdot) & 
\stackrel{\mbox{\tiny{$\hat \partial_0$}}}{\longrightarrow} &
 \hat \pi_{k+1} (X_c,  \centerdot),
 \end{array}
 \end{equation}

 With these notations we have the following:

\begin{proposition}\label{p:2}
Let $(\bX_{D_i}, X_{c_i})$ be $(k+1)$-connected, $k\ge 1$,  for all $i$.  
Then $\im \hat \partial_0 = \im (\sum_i \hat \partial_i)$ in 
diagram \rm (\ref{eq:3})\it.
 \end{proposition}
 \begin{proof}
We use Hurewicz maps between the relative homotopy and homology groups (denoted $\cH_i$ and $\cH$ below). We have the following commutative diagram:
\[ \begin{array}{ccc}
 \oplus_i \hat\pi_{k+2} (\bX_{D_i}, X_{c_i}, \centerdot)  & \stackrel{\sum 
 \hat{\iota_i}_\sharp}{\longrightarrow} & \hat\pi_{k+2} (\bX_{\cD},X_c, \centerdot)  \\
 \mbox{\tiny{$\oplus_i \cH_i$}} \downarrow  & \ &   \downarrow 
 \mbox{\tiny{$\cH$}}  \\
\oplus_i H_{k+2} (\bX_{D_i}, X_{c_i})  & \stackrel{\sum   {\iota_i}_*}{\longrightarrow} & H_{k+2} (\bX_{\cD},X_c). 
 \end{array}
 \]
The additive notations are used because the relative homotopy groups 
are abelian for $k\ge 1$.

 Under our hypothesis and by the relative Hurewicz isomorphism theorem (see e.g. \cite[p. 397]{Sp}), we get that
 the Hurewicz map $\cH_i : \hat\pi_{k+2} (\bX_{D_i}, X_{c_i}, \centerdot)\to
 H_{k+2} (\bX_{D_i}, X_{c_i})$ is an isomorphism, for any $i$. The same is
 $\cH$, since $(\bX_\cD, X_c)$ is $(k+1)$-connected, by
Lemma \ref{l:prop}.

 Next, $\sum {\iota_i}_*$ is a homology excision, hence an isomorphism too. It
 follows that
 $\sum \hat{\iota_i}_\sharp$ is an isomorphism too, which proves our claim.
\end{proof}

 To complete the proof of our main theorem, we need the following:

\smallskip
 \noindent
$\bullet$  \ to prove that $(\bX_{D_i}, X_{c_i})$ is $(k+1)$-connected, 
$k\ge 0$, for all
$i$. 

\noindent 
$\bullet$ \  to find the image of the map $\hat \partial_i : \hat\pi_{k+2} (\bX_{D_i}, X_{c_i}, \centerdot) \to \hat 
\pi_{k+1} (X_{c_i}, \centerdot)$, $k\ge 1$, for all
$i$.

\smallskip

First thing we do is to show that we may replace $\bX_{D_i}$ by $\bX^*_{D_i}:= 
\bX_{D_i}\m \Sing_\cS p$, using the homotopy depth assumption, the definition
of which we recall here.
\begin{definition}\label{d:hd}
For a discrete subset $\Phi \subset \bX$, one sais that the {\bf homotopy
 depth} of $\bX$ at $\Phi$, denoted by $\hd_\Phi \bX$, is greater or equal to $q+1$ if, at any 
point $\alpha\in \Phi$, 
there is an arbitrarily small neighbourhood
$\cN$ of $\alpha$ such that the pair $(\cN, \cN\m \{\alpha\})$ is $q$-connected.
\end{definition}
\begin{lemma}\label{l:2}
If $\hd_{\bX \cap\Sing_\cS p} \bX \ge q+1$, for $q\ge 1$, then the inclusion of pairs 
$(\bX^*_{D_i}, X_{c_i}) \stackrel{j}{\hookrightarrow} (\bX_{D_i}, X_{c_i})$ is a 
$q$-equivalence, for all $i$.
\end{lemma}
 \begin{proof}
 In the exact sequence of the triple $(\bX_{D_i},\bX^*_{D_i}, X_{c_i})$, it is 
sufficient to prove that $(\bX_{D_i},\bX^*_{D_i})$ is $q$-connected, for all
$i$. By the homotopy depth assumption, for all $j$, the pair $(\bX_{D_i}\cap
B_{ij}, \bX_{D_i}\cap B_{ij}\m \{a_{ij})$ is $q$-connected. 
By Switzer's result for CW-complexes  \cite[6.13]{Sw}, the first space is
obtained from the second by attaching cells of dimensions $\ge q+1$.
It follows that $\bX_{D_i}$ is obtained from $\bX^*_{D_i}$
by replacing each $B_{ij}\m \{a_{ij}$ with $B_{ij}$, which amounts, as we have
seen, to attaching cells of dimensions $\ge q+1$. So our claim is proved.
 \end{proof}  


\begin{corollary}\label{c:2}
If $\hd_{\bX \cap\Sing_\cS p} \bX \ge k+3$, $k\ge 0$, then, for all $i$:
\[ \im (\hat\partial_i : \hat\pi_{k+2}(\bX_{D_i}, X_{c_i}, \centerdot) \to \hat\pi_{k+1}(X_{c_i},  \centerdot)) = \im 
(\hat\partial'_i : \hat\pi_{k+2}(\bX^*_{D_i}, X_{c_i}, \centerdot) \to \hat\pi_{k+1}(X_{c_i},  \centerdot)).
\]
\end{corollary}
\begin{proof}  By Lemma \ref{l:2}, $j_\sharp : 
\pi_{k+2}(\bX^*_{D_i}, X_{c_i}, \centerdot) \to \pi_{k+2}(\bX_{D_i}, X_{c_i},
\centerdot)$ is surjective. Since $\partial'_i = \partial_i \circ j_\sharp $,
we get our claim after taking the quotient by the $\pi_1(X_{c_i},\centerdot)$.
\end{proof} 

 The last step in the proof of Theorem \ref{t:main} consists of the following
 two key results.

\begin{proposition}\label{p:3}
If $(X_{c_i}, X^*_{a_i})$ is $k$-connected, $k\ge 0$, then
 $(\bX^*_{D_i}, X_{c_i})$ is $(k+1)$-connected.
\end{proposition}
\begin{proof} 
 
 We shall suppres the lower indices $i$ in the following.

\noindent
 Let $D^* := D\m \{ a\}$. We retract $D^*$ to the circle $\partial D$ and cover 
this circle by two arcs $I\cup J$. We have the homotopy equivalences
 $\bX_{D^*} \stackrel{\h}{\simeq} \bX_I \cup \bX_J$ 
and 
$X_c \stackrel{\h}{\simeq} \bX_J$. We remind that $\bX_D \m \cup_j B_{ij}$ is 
the total space of the trivial fibration 
 $p_| \colon \bX_D \m \cup_j B_{ij}\to D$, its fibre being homotopy equivalent
 to $X_c \m \cup_j B_{ij}$ (which we have also denoted by $X^*_a$).
Then the following inclusions are homotopy equivalences of pairs and they
induce isomorphisms in all relative homotopy groups:

\[ \begin{array}{c}
(\bX^*_D , X_c) \hookleftarrow (\bX_{D^*} \cup (\bX_D \m \cup_j B_{ij}), X_c)
\hookrightarrow
(\bX_{D^*} \cup (\bX_D \m \cup_j B_{ij}), \bX_J\cup (\bX_D 
\m \cup_j B_{ij})) \hookleftarrow \\
\hookleftarrow (\bX_I \cup \bX_J \cup (\bX_D \m \cup_j B_{ij}), \bX_J \cup 
(\bX_D \m \cup_j B_{ij})).\end{array}\]
 
Let us consider the following homotopy excision:
\begin{equation}\label{eq:htex}
(\bX_I \cup \bX_J \cup (\bX_D \m \cup_j B_{ij}), \bX_J \cup 
(\bX_D \m \cup_j B_{ij})) \hookleftarrow (\bX_I, \bX_{ 
\partial I} \cup (\bX_I \m \cup_j B_{ij})).
\end{equation} 

The right hand side pair is homotopy equivalent to 
$(\bX_I, X_c \times 
\partial I \cup X^*_a \times I)$, which is
the product of pairs $(X_c, X^*_a) \times (I, \partial I)$.
Since $(X_c, X^*_a)$ is $k$-connected by our hypothesis, this product is 
$(k+1)$-connected\footnote{This follows from the fact that the spaces are
  CW-complexes and by using Switzer's \cite[Proposition 6.13]{Sw}.}.
 
 By the Blakers-Massey theorem (cf Gray \cite[Cor; 16.27]{Gr}), it will follow that the excision 
(\ref{eq:htex}) is a $(k+1 +q)$-equivalence,
 provided that the pair $(\bX_J \cup (\bX_D \m \cup_j B_{ij}), \bX_{\partial
   I} \cup (\bX_I \m \cup_j B_{ij}))$ is $q$-connected. If $q\ge 0$, then this
 proves our claim that the original pair $(\bX^*_D , X_c)$ is  $(k+1)$-connected.

It therefore remains to evaluate the level $q$. Since the inclusions $\bX_I\m
\cup_j B_{ij} \hookrightarrow \bX_D \m \cup_j B_{ij}$ and $X_c\hookrightarrow \bX_J\m
\cup_j B_{ij} \hookrightarrow \bX_D \m \cup_j B_{ij}$ are homotopy
equivalences, it follows that 
the following inclusions of pairs are homotopy 
equivalences and therefore induce isomorphisms in all relative homotopy groups: 

\[ \begin{array}{c} (\bX_J \cup (\bX_D \m \cup_j B_{ij}), \bX_{\partial I} \cup (\bX_I \m 
\cup_j B_{ij})) 
\hookrightarrow 
(\bX_J \cup (\bX_D \m \cup_j B_{ij}), \bX_{\partial I} \cup (\bX_D \m \cup_j B_{ij})) 
\hookleftarrow \\
(\bX_J \cup (\bX_J \m \cup_j B_{ij}), \bX_{\partial J} \cup (\bX_J \m \cup_j B_{ij})) 
\hookleftarrow 
(\bX_J, \bX_{\partial J} \cup (\bX_J \m \cup_j B_{ij})). 
\end{array} \]
We also have:
\[ (\bX_J, \bX_{\partial J} \cup (\bX_J \m \cup_j B_{ij})) \stackrel{\h}{\simeq} \\
(\bX_c \times J, X_c \times \partial J \cup X^*_a \times J)
= (X_c, X^*_a) 
\times (J, \partial J).\]
Since  $(X_c, X^*_a)$ is supposed $k$-connected, this implies that $q$ is at
least $k+1$.
\end{proof}

\begin{proposition}\label{p:4}
If $(X_{c_i}, X^*_{a_i})$ is $k$-connected, $k\ge 2$, then \\
$\im \hat\partial'_i = \im (\widehat{\hvar_i} : \hat\pi_{k+1}(X_{c_i}, X^*_{a_i}, 
\centerdot) \to \hat\pi_{k+1}(X_{c_i}, \centerdot))$.
\end{proposition}
\begin{proof} 
 
 As in the above proof, we suppres the lower indices $i$.
The final purpose is to show the commutativity of the
diagram (\ref{eq:varcomm}) below. The idea is to compare it with the analogous
diagram for homology
groups, via the Hurewicz maps. However, we cannot do it directly, since
we do not have the $k$-connectivity of $X_c$ from the upper right corner 
of (\ref{eq:varcomm}).
 We only have that the pair $(X_{c}, X^*_a)$ is $k$-connected. So we give
a proof in two steps. 
 
 Let us consider the
following diagram:
\begin{equation}\label{eq:ht}
 \begin{array}{lcc}
\pi_{k+2} (\bX^*_{D} , X_{c}\times I, \centerdot)\simeq 
\pi_{k+2} (\bX^*_{D} , X_{c}, \centerdot) & \stackrel{\partial''}{\longrightarrow} & \pi_{k+1} 
(X_{c}, X^*_a, \centerdot) \\
\mbox{\tiny exc} \uparrow  & \  & \uparrow \mbox{\tiny $\mu$} 
\\
 \pi_{k+2} (\bX_I, X_{c} \times \partial I \cup X^*_{a} \times I,
 \centerdot) & \stackrel{\simeq}{\leftarrow} & \pi_{k+1} 
(X_{c},X^*_{a}, \centerdot),
\end{array} 
\end{equation}
where $\mu$ is the morphism defined in (\ref{eq:mu}) and $\partial''$ denotes
the boundary morphism of the homotopy exact sequence of the triple 
$(\bX^*_{D} , X_{c}, X^*_{a})$. The fundamental group $\pi_1(X^*_a,
\centerdot)$ acts on the groups at the right hand side of this diagram
and via isomorphisms on the groups at the left  hand side. Indeed we have the
isomorphisms $\pi_1(X^*_a,\centerdot) \simeq \pi_1(X_c,\centerdot) 
\simeq \pi_1(X_c\times \partial I \cup X^*_{a} \times I,\centerdot)$ which
follow from the connectivity hypothesis on  $(X_{c}, X^*_a)$ and from the
inclusion of spaces $X^*_{a} \times I\subset X_{c} \times \partial I \cup
X^*_{a} \times I \subset X_{c}\times I$.
Therefore the natural diagram (\ref{eq:ht}) passes to the quotient by the $\pi_1$-action 
and yields:
\begin{equation}\label{eq:mucomm}
 \begin{array}{lcc}
\hat\pi_{k+2} (\bX^*_{D} , X_{c}\times I, \centerdot)\simeq 
\hat\pi_{k+2} (\bX^*_{D} , X_{c}, \centerdot) & \stackrel{\hat\partial''}{\longrightarrow} & \hat\pi_{k+1} 
(X_{c}, X^*_a, \centerdot) \\
\widehat{\mbox{\tiny exc}} \uparrow  & \  & \uparrow \mbox{\tiny $\hat\mu$} 
\\
 \hat\pi_{k+2} (\bX_I, X_{c} \times \partial I \cup X^*_{a} \times I,
 \centerdot) & \stackrel{\simeq}{\leftarrow} & \hat\pi_{k+1} 
(X_{c},X^*_{a}, \centerdot),
\end{array} 
\end{equation}

We claim that (\ref{eq:mucomm}) is a commutative diagram.
In homology, there is the following diagram: 
\begin{equation}\label{eq:homol}
\begin{array}{ccc}
H_{k+2} (\bX^*_{D} , X_{c}) & \stackrel{\partial''}{\longrightarrow} & H_{k+1} 
(X_{c}, X^*_a) \\
\mbox{\tiny exc} \uparrow  & \  & \uparrow \mbox{\tiny $h_* - \id$} 
\\
 H_{k+2} (\bX_I, X_{c} \times \partial I \cup X^*_{a} \times I) & \stackrel{\simeq}{\leftarrow} & H_{k+1} 
(X_{c},X^*_{a}),
\end{array} 
\end{equation} 
  which is a Wang type diagram, and therefore it is {\em commutative}. 
Milnor explained a simpler version of it in
  \cite[pag. 67]{Mi} and his explanation holds true in the relative homology. 
   We have used in \cite[\S 3]{Ti-xiv} an  equivalent
  commutative diagram,  which differs from (\ref{eq:homol}) by the upper
  right corner. 

 Let's now prove the claim. 
We parallel the homotopy groups diagram (\ref{eq:mucomm}) by the homology
groups diagram (\ref{eq:homol}) and connect them by the Hurewicz
homomorphism. We get in this way a ``cubic'' diagram. 
Notice that the homology version of the homotopy map $\mu$ is just $h_* -\id$,
 by the very definition of $\mu$ given in \S \ref{ss:var}.

By the naturality of the Hurewicz morphism, all
the maps between homotopy groups are in commuting diagrams with their homology 
versions (as ``faces'' of the cubic diagram). Moreover, the Hurewicz theorem
may be applied each time. The connectivity conditions that the 4 Hurewicz
morphisms become isomorphisms are fullfilled: $(X_{c}, X^*_{a})$ is
$k$-connected by hypothesis and  $(\bX^*_{D}, X_{c})$ is $(k+1)$-connected by 
Proposition \ref{p:3}.
 This allows one to identify the 
homotopy diagram (\ref{eq:mucomm}) to the 
homology one (\ref{eq:homol}) and since the later commutes, it 
is the same with the former.  Our claim is therefore proved.

\smallskip 
 We now claim that the following
 diagram, which differs from (\ref{eq:mucomm}) by the upper right corner only,
 is commutative too:
\begin{equation}\label{eq:varcomm}
 \begin{array}{ccc}
\hat\pi_{k+2} (\bX^*_{D} , X_{c}, \centerdot) & \stackrel{\hat\partial'}{\longrightarrow} & \hat\pi_{k+1} 
(X_{c}, \centerdot) \\
\widehat{\mbox{\tiny exc}} \uparrow & \  & \uparrow \mbox{\tiny $\widehat{\hvar}$} 
\\
 \hat\pi_{k+2} (\bX_I, X_{c} \times \partial I \cup X^*_{a} \times I, \centerdot) & \stackrel{\simeq}{\leftarrow} & \hat\pi_{k+1} 
(X_{c},X^*_{a}, \centerdot),
\end{array} 
\end{equation}

Let us justify why we may replace the morphisms $\mu$ and $\partial''$ in 
(\ref{eq:mucomm}) by the morphisms $\hvar$ and $\partial'$ respectively. We
 have the following diagram:

\begin{equation}\label{eq:pullback}
 \begin{array}{c}
 
 \xymatrix{
    & &  \hat\pi_{k+1}( X_c, X^*_a, \centerdot) \\
        \hat\pi_{k+2}( \bX^*_D, X_c, \centerdot) \ar@/^/[urr]^{\hat\partial''}
     \ar[r]^{\hat\partial'} & \hat\pi_{k+1}( X_c,  \centerdot) \ar@{>}[ur]|-{\hat{j_\sharp}} &   \\
     & \hat\pi_{k+1}( X_c, X^*_a, \centerdot) \ar@/_/[uur]_{\hat\mu} \ar[u]_{\widehat{\hvar}} &   }    
 \end{array}
\end{equation}
     
 The two triangles of this diagram are commutative since: 1). the right hand
 side triangle coincides with diagram 
(\ref{eq:diag}) after having taken quotients by the action of
 $\pi_1(X_c,  \centerdot)$, and 2).  the upper-left triangle is commutative by
 the naturality of the boundary morphism.    
 
 Let now $\alpha \in \pi_{k+1}(X_{c},X^*_{a}, \centerdot)$.
 As explained in \S \ref{ss:var}, an element of the form $\mu(\alpha)$ is
 homotopic, relative to $X_a^*$, to a certain element of the absolute group 
$\pi_{k+1}(X_{c}, \centerdot)$, which we have denoted by $\hvar(\alpha)$.
  Furthermore, let us remark that, by the commutativity of diagram (\ref{eq:mucomm}), $\hat\mu(\alpha) = 
\hat\partial''(\beta) \in \hat\pi_{k+1}(X_{c}, X^*_{a},
 \centerdot)$ for some $\beta\in \hat\pi_{k+2}(\bX_D^*, X_{c}, \centerdot)$.
We get  $\widehat{\hvar}(\alpha) = 
\hat\partial'(\beta)$.
  
 We have thus shown that one can ``pull-back'' (as in the figure
 (\ref{eq:pullback})) the upper right corner of the
 diagram (\ref{eq:mucomm}).

 The proof of our statement  follows now from the commutativity of the diagram
 (\ref{eq:varcomm}), and the fact that the botom row and the excision from the
 left are both isomorphisms.

\end{proof}

We may now conclude the proof of Theorem \ref{t:main} as follows:

\smallskip
\noindent
$\bullet$ The claim (a) follows by chaining together
Lemmas \ref{l:im}, \ref{l:1}, \ref{l:prop}, \ref{l:2} (for $q=k+1$), 
 and finally Proposition \ref{p:3}. 

\noindent
$\bullet$ 
The claim (b) follows by chaining together Lemma \ref{l:im}, Proposition \ref{p:2}, 
Corollary \ref{c:2} and Proposition \ref{p:4}.

\section{Some special cases of the main theorem}\label{s:further}

The Lefschetz method for generic pencils of hyperplanes
and the comparison with the axis of the pencil
have been used before by several authors  in order to prove 
connectivity results of increasing generality.  
Our Theorem \ref{t:main}(a) improves \cite{Ti-lef} by weakening
 the hypotheses. As we have explained in {\em loc.cit.}, our result \cite{Ti-lef}
 recovered in its turn other previous connectivity results
\cite{La,Ch-ens, Ch, Ey}.

  In the following we show
what becomes Theorem \ref{t:main}(b) in some particular
cases. 
 
\subsection{Aspects of the homotopy depth condition}\label{ss:depth}
The condition (iii)  of Theorem \ref{t:main} is satisfied 
as soon as $X$ is a {\em locally complete intersection} of dimension $\ge
k+2$.  This follows from the fact that the 
{\em rectified homotopical depth}\footnote{see \cite{Ti-lef} and the 
  source \cite{HL} for the definition and properties of this notion which goes
back to Grothendieck; in particular, $\rhd$ does not depend on the
stratification.} is $\ge k+2$ for such an
 $X$ (cf \cite{HL}).  
 Indeed, one can easily show that, more generally, the condition $\rhd \bX \ge k+2$
implies the condition (iii), i.e. that  $\hd_{\bX \cap \Sing_\cS p} 
\bX \ge k+2$.\footnote{see also \cite{Ti-xiv} for a similar remark in
  homology.}
Indeed, $\rhd X \ge k+2$ implies $\rhd \bX \ge k+2$, 
since we may apply  \cite[Theorem 3.2.1]{HL} to the hypersurface $\bX$ of $X\times 
\bP^1$. This in turn implies 
$\hd_v \bX \ge k+2$, for any point $v\in \bX$, by definition. 

 Moreover,  $\rhd \bX \ge k+2$ implies that the pair $(\bX_{D_i}, X_{c_i})$ is
 $(k+1)$-connected, by 
\cite[Proposition 4.1]{Ti-lef}. This provides a shortcut to the proof of 
Part (a) of Theorem \ref{t:main}, in the sense that the condition (ii) is not
needed, and therefore Lemma \ref{l:2} and Proposition \ref{p:3} are useless.

In what concerns Part (b) of Theorem  \ref{t:main}, the condition 
  $\rhd X \ge k+3$ implies not only (iii'), as shown above, but also 
the condition (ii), which is important in the proof of Part (b). 
Indeed, the $k$-connectivity of the pair $(X_{c_i}, X^*_{a_i})$ follows from 
a sequence of homotopy excisions, which reduces the problem to proving
  the $k$-connectivity of each pair $(X_{c_i}\cap B_{ij}, X_{c_i}\cap \partial
  B_{ij})$. This is true since $X_{c_i}\cap B_{ij}$ is $(k+1)$-connected
and $X_{c_i}\cap \partial B_{ij}$ is $k$-connected. The former follows from 
\cite[Proposition 3.4]{Ti-lef} which is based on \cite[Corollary 4.2.2]{HL}.
The proof goes as follows:  $\rhd X \ge k+3$ implies $\rhd X_{a_i} \ge k+2$, by  
\cite[Theorem 3.2.1]{HL}. This implies $\hd_{a_{ij}} X_{a_i} \ge k+2$,
which  shows that $X_{a_i}\cap \partial B_{ij}$ is $k$-connected. 
Finally, we already know that 
$X_{a_i}\cap \partial B_{ij}$ is homeomorphic to $X_{c_i}\cap \partial
  B_{ij}$, so this ends the proof of our claim.

\subsection{Pencils on $X$ having no singularities in the axis}\label{ss:notaxis}

Let us consider the situation $(A\times \bP^1)\cap\bX \cap \Sing_\cS p
=\emptyset$. We say that ``$X$ does not contain singularities along the axis of
the pencil''. Even so,
 there might still be singularities in the axis on $Y$, and they certainly
 influence the topology of the pencil on $X$.  However, this situation  be 
further investigated as follows, see also \S \ref{ss:poly} for an interesting case.
 We claim that the
 homotopy depth condition (iii), resp. (iii'), may be replaced by
the following more general, global condition:
\begin{equation}\label{eq:glob}
(X, X\m \Sigma) \mbox{ is $(k+2)$-connected, resp. $(k+3)$-connected, } 
\end{equation}
where $\Sigma :=  \sigma(\Sing_\cS p) \subset Y$.
Then our  Theorem \ref{t:main}(a) specializes to a connectivity statement
 which recovers Eyral's main result \cite{Ey}, where such a condition was used
 for generic pencils.
  Indeed, we notice that the homotopy depth condition was used only when
  comparing $\bX_{D_i}$ to $\bX^*_{D_i}$. We may therefore get rid of this
  comparison by replacing everywhere
 in the proof of Theorem \ref{t:main} the
 space $X$ by the space $X\m \Sigma$. In particular, this reduces to
 tautologies  Lemma \ref{l:2} and Corollary \ref{c:2}. We then get the
 conclusion of Theorem \ref{t:main} for the inclusion
$X_c \hookrightarrow X\m\Sigma$ instead of the inclusion $X_c \hookrightarrow X$. 
But, at this stage, we may just
plug in $X$ in the place of $X\m\Sigma$ since the condition (\ref{eq:glob}) tells that
$X\m\Sigma$ and $X$ have isomorphic homotopy groups up to $\pi_{k+1}$, and
this is all we need. 

Still in case $(A\times \bP^1)\cap\bX \cap \Sing_\cS p =\emptyset$, let us
observe that if the condition (i) of Theorem \ref{t:main} is fulfilled, then
the condition (ii) becomes equivalent to:
\begin{equation}\label{eq:spec}
(X^*_{a_i}, A\cap X^*_{a_i}) \mbox{ is $(k-1)$-connected}.
\end{equation}

Indeed, when there are no singularities in the axis, we have $A\cap X^*_{a_i} = A\cap X_c$, 
for any $i$ and our claim follows from the exact sequence of the triple 
$(X_c, X^*_{a_i}, A\cap X^*_{a_i})$.

We have to point out that even if there are ``no singularities 

\subsection{The classical case: generic pencils on the projective space}\label{ss:classic}
  From the preceding observations on pencils without singularities in the
 axis, one may  derive the following statement for {\em complements of arbitrarily singular
 subspaces} in $\bP^n$.
   Its proof goes by induction, based on the abundance of generic pencils,
 i.e. pencils  with no singularities along the axis and such that the axis $A$ is
 transversal to all strata and is not included into the subspace $V$. 
The proof actually follows the pattern of the homology proof discussed in \cite[\S
 4.3]{Ti-xiv}, supplemented by our homotopy considerations in \S
 \ref{proof} above. 
\begin{corollary}\label{c:main}
Let $\cY= \bP^n$. Let $V\subset \cY$ be a singular complex algebraic 
subspace, not necessarily irreducible.
For any hyperplane $H\subset \cY$ transversal to all strata of $V$, we have:
 \begin{enumerate}
\rm \item \it  $H\cap (\cY \m V) \hookrightarrow \cY \m V$ is a $(n+\codim V -2)$-equivalence.
\rm \item \it There exist generic pencils of hyperplanes having $H$ as a
generic member. 
If $\dim V\le 2n-5$ then  the kernel  of the surjection $\pi_{n+\codim V -2}(H\cap 
(\cY\m V), \centerdot) \twoheadrightarrow \pi_{n+\codim V -2} (\cY \m V,
\centerdot)$ is 
generated by the images of the variation maps $\hvar_i$ of such a generic pencil.  
\fin \end{enumerate}
 \end{corollary}
\begin{proof}
The part (a) of this Corollary is well-known, 
see e.g. \cite{Ch-ens,Ey}, and its proof goes by  induction on $\dim V$, since
by repeatedly slicing we arrive to the case $\dim V=0$. 
The conditions of Theorem \ref{t:main} are satisfied at each induction step, as follows.
 The generic slice $X_c$ is path connected at each step.
Condition (iii), resp. (iii'), is empty at 
every step, since a generic pencil has no singularities 
in $\bP^{n-j} \m V$. Condition (ii) is implied by (i) since it is equivalent to
(\ref{eq:spec}), as shown above. The condition (i) in case $\dim V=0$,  is
 the level of  connectivity of the pair $(\bP^{n-\dim V -1}, \bP^{n-\dim V-2})$, 
which is known to be equal to $k= 2(n-\dim V -1)-1$.
Then Theorem \ref{t:main}(a) gives the level of connectivity $k+1$
of the higher pair $(\bP^{n-\dim V} \setminus V\cap \bP^{n-\dim V},
\bP^{n-\dim V-1})$, which becomes the condition (i) for the
next induction step. This proves part (a). Part (b) also follows from  Theorem
\ref{t:main}(b),
since we have seen that the conditions are inductively fulfilled. In addition,
we have
to insure that the initial step verifies the requirement $k\ge 2$.
This amounts to $\dim V \le 2n-5$. 
\end{proof}
 The part (b) of this Corollary seems to be a new result since we consider the general
 situation of complements of arbitrarily singular subspaces. There is 
 a single case which is not covered by the condition  $\dim V \le 2n-5$,
 namely the case $n=3$,  $\dim V =2$ (since $\hvar$ is not well defined for $\pi_2$). 
 Of course, the case $n=2$ does not
 matter here since it is covered by the classical Zariski-van Kampen theorem.
 
 In the particular case $V\subset \bC^n$ is a hypersurface
 with at most isolated singularities and transversal to the hyperplane at infinity, a similar type
of result, but with different background,  has been proved by Libgober
 \cite{Li}, including at the $\pi_2$ level.


\subsection{A complementary case}\label{ss:poly}
It appears that the case when $V$ contains a member of the pencil, 
which was excluded in our main theorem (since
 $A\subset V$ in this case), can also be
 treated. For homology groups, this was explained in  \cite[\S
 4-5]{Ti-xiv}. Our main
 result can be reformulated, with much less restrictive conditions.
 
\begin{theorem}\label{t:2} 
  Let $h :Y\dashrightarrow \bP^1$ define a pencil with isolated
  singularities, such that $V$ 
contains a member of the pencil. Le $X_c$ be path connected and let
$(X_c, X^*_{a_i})$ be  $k$-connected, for a general member $X_c$ and any atypical one
  $X^*_{a_i}$.  Then:
  
\begin{enumerate}
\rm \item \it If $k \ge 0$ and $(X, X\m \Sigma)$ is $(k+1)$-connected, then the inclusion 
$X_c \hookrightarrow X$ is a $(k+1)$-equivalence.
\rm \item \it If $k\ge 2$ and if $(X, X\m \Sigma)$ is $(k+2)$-connected, then
the kernel of the surjection $\pi_{k+1}(X_c) \twoheadrightarrow \pi_{k+1}(X)$
is generated by the images of the variation maps $\hvar_i$. 
\end{enumerate} 
\fin
\end{theorem}
\begin{proof}
The proof follows by revisiting the one of Theorem \ref{t:main}. We then
 observe that in our case
$X_\cD$ is just homotopy equivalent to $X$, according to the definition of
 $\cD$. We have used the condition (i) only to
compare these two spaces, in Lemma \ref{l:im}, so condition (i) does not occur
 any more in our situation.  
Let us further  remark that the condition $(A\times \bP^1)\cap\bX \cap \Sing_\cS p
=\emptyset$  considered in \S \ref{ss:notaxis} is satisfied here, and so we may 
use condition (\ref{eq:glob}) istead of (iii) and (iii'), respectively.
\end{proof}

\end{document}